\documentclass[amsfonts,amsmath,tabularx,11pt]{article}
\newtheorem{theorem}{Theorem}[section]
\newtheorem{lemma}[theorem]{Lemma}
\newtheorem{result}[theorem]{Result}
\newtheorem{proposition}[theorem]{Proposition}
\newtheorem{corollary}[theorem]{Corollary}
\newtheorem{definition}[theorem]{Definition}
\newenvironment{definition-new}{\begin{definition} \em}{\end{definition}}
\newtheorem{remark}[theorem]{Remark}
\newenvironment{remark-new}{\begin{remark} \em}{\end{remark}}
\newtheorem{example}[theorem]{Example}
\newenvironment{example-new}{\begin{example} \em}{\end{example}}

\newenvironment{notation-new}{\begin{remark} \em}{\end{remark}}

\newenvironment{agr-new}{\begin{remark} \em}{\end{remark}}

\makeatletter \@addtoreset{equation}{section} \makeatother

\makeatletter \@addtoreset{figure}{section} \makeatother
\usepackage{amsmath}
\usepackage{amsfonts}

\begin{document}
\begin{center}
{\textbf {\Large Harmonicity of the complex structure on product of trans-Sasakian manifolds }}\\[2pt]
 {\bf Nidhi Yadav}\footnote{$^{,3}$
Department of Mathematics \& Statistics, Dr. Harisingh Gour Vishwavidyalaya,
Sagar-470 003, M.P. INDIA \newline
Email: nidhiyadav.bina@gmail.com, rkgangele23@gmail.com}, {\bf Punam Gupta}\footnote{
School of Mathematics, Devi Ahilya Vishwavidyalaya, Indore-452 001, M.P. 
INDIA\newline
Email: punam2101@gmail.com} and {\bf R.K. Gangele}$^3$

\end{center}

\noindent
{\bf Abstract:} In this paper, we investigate the transverse geometry of trans-Sasakian manifolds and present several significant findings. We analyze the Levi-Civita connection associated with the metric on the product manifold of two trans-Sasakian manifolds. We outline the conditions under which the complex structure is harmonic on the product manifold. Notably, we also offer valuable insights into the harmonicity of the complex structure within the context of astheno-K\"ahler manifolds.
\vskip.2cm
\noindent {\bf MSC:} 53C05, 53C15, 53C25, 53D15. 
\vskip.2cm
\noindent
{\bf   Keywords:} Trans-Sasakian manifold, Harmonic map, Astheno-K\"ahler manifold 

\section{Introduction}
K\"ahler manifolds play central roles in complex and algebraic geometry, Hodge theory, and in theoretical physics, such as quantum mechanics and supersymmetry.

A compact complex manifold \( M^n \) that possesses a K\"ahler metric \( g \) must exhibit a particularly unique topological structure. According to the classical findings of Deligne et al. \cite{del}, the rational homotopy type of such a manifold is dictated by its cohomology ring. Additionally, the Hodge and Lefschetz decomposition theorems impose significant constraints on the Betti numbers of \( M^n \). For instance, all odd Betti numbers \( b_{2k-1} \) must be even, while all even Betti numbers \( b_{2k} \) must be positive, applicable to any degree from 0 up to the real dimension of \( M^n \). Consequently, many compact complex manifolds are unable to support a K\"ahler metric.

Over the last thirty years, astheno-K\"ahler manifolds have attracted significant interest from geometers worldwide, highlighting their potential importance in both differential and algebraic contexts. The term ``astheno-K\"ahler manifold," which means weak in Greek, was first introduced by Jost and Yau \cite{jost,jost2} in their seminal 1993 paper. A crucial contribution by Matsuo and Takahashi \cite{km} established that every compact balanced astheno-K\"ahler manifold is K\"ahler, indicating that compact Hermitian-flat astheno-K\"ahler manifolds also fall into this category. Their research yielded several key results, including that any product manifold formed from curves and surfaces is astheno-K\"ahler, and that the product of two compact Sasakian manifolds in three dimensions, equipped with a specific Hermitian structure and K\"ahler form, is also astheno-K\"ahler. Furthermore, if one manifold is a three-dimensional compact Sasakian manifold and the other is a compact cosymplectic manifold of dimension at least three, their product retains the astheno-K\"ahler property. In another study, Matsuo \cite{kmm} identified astheno-K\"ahler metrics on Calabi-Eckmann manifolds, which are principal \( T^2 \)-bundles over \({ \mathbb CP}^n \times {\mathbb CP}^m\), thus enabling the development of astheno-K\"ahler structures on torus bundles over K\"ahler manifolds. More on astheno-K\"ahler manifolds can be observed in \cite{Gupta}.

To grasp the complexities of non-K\"ahler Hermitian geometry, various methodologies have been employed. These diverse approaches provide insights that enhance our understanding of this specialized field, allowing for a more comprehensive exploration of its unique properties and implications.

Out of them, one alternative method involves identifying an optimal almost complex structure on a given Riemannian manifold $(M, g)$ from a variational perspective. Assuming the existence of at least one almost complex structure that is compatible with the metric, the approach focuses on functionals defined over the set of orthogonal almost complex structures on $M$, seeking to find their extrema. A prominent functional \cite{wood} that has been studied in previous works is the Dirichlet energy functional $E$, expressed as 
\begin{align}
E(J):=\int \vert|\nabla J\vert|^{2}vol_g,
\label{eq-1a}
\end{align} 
when $M$ is compact. The initial step in locating (local) minima entails calculating the critical points of the functional $E$, referred to as harmonic almost complex structures. It has been established that an orthogonal almost complex structure $J$ is harmonic if and only if the equation $[J, \nabla^*\nabla J] = 0$ holds, where $\nabla^*\nabla J$ denotes the rough Laplacian of $J$, defined as $\nabla^*\nabla J = Tr\nabla^2 J$. In the case of a non-compact manifold, one can either adopt the Euler-Lagrange equation $[J, \nabla^*\nabla J] = 0$ as the definition of harmonicity or define the energy on an open subset with compact closure and consider variations with compact support within this subset, leading to an equation that mirrors the compact case.

It is evident from equation (\ref{eq-1a}) that K\"ahler structures possess harmonic properties, functioning as absolute minimizers. In a broader Hermitian context, research presented in \cite{gan} has demonstrated that the complex structure of a balanced Hermitian manifold is also harmonic, as is the case for LCK manifolds when the complex dimension exceeds two. Additionally, examples of harmonic almost complex structures can be found on Calabi-Eckmann manifolds that are endowed with the product of round metrics \cite{wood}. More on harmonic almost complex structures can be observed in various studies, such as those referenced in \cite{bor,lo}.

In our recent study \cite{nidhi}, we investigate the conditions
under which astheno-K\"ahler structures can be identified on the product of two compact trans-Sasakian manifolds of dimensions greater than $2$. Inspired by these studies, this work investigates the instances in which astheno-K\"ahler manifolds possess a harmonic complex structure relative to their metric.

The article is structured as follows. In section $2$, we recall basic concepts  and results on almost complex and almost contact
manifolds. In section $3$, we study the transverse geometry of trans-Sasakian manifolds and present some useful results. In section $4$,
 we study the Levi-Civita connection of the metric $g_{a,b}$ on the product $M_1\times M_2$,where $M_1$
 and $M_2$ are trans-Sasakian manifolds. In the last section, we give the conditions when $J_{a,b}$ is harmonic on
 $(M_1\times M_2,g_{a,b})$. In particular, we provide some useful results on the harmonicity of complex structure on the astheno-K\"ahler manifold.
\section{Preliminaries}
In this section, we review some basic concepts and results on almost complex and almost contact manifolds.
\subsection{Almost complex manifolds}
An almost complex manifold $(M^{2n},J)$ \cite{ballmann} together with a compatible
Riemannian metric $g$, that is, $g(JX, JY ) = g( X, Y)$
for all vector fields $X$, $Y$ on $M$, is called an almost Hermitian manifold and metric $g$ is called an almost Hermitian metric. A complex manifold $M$ \cite{ballmann} together with a compatible
Riemannian metric $g$ is called a Hermitian manifold and metric $g$ is called a Hermitian metric. 
The alternating $2$-form
 $$\Omega(X, Y ) = g(JX, Y)$$
is called the associated K\"{a}hler form. We can retrieve $g$ from $\Omega$
$$g( X, Y ) = \Omega(X, JY ).$$
\noindent If $\Omega $ is closed, then $(M,g)$ is known as a K\"{a}hler
manifold \cite{ballmann} and $g$ is a K\"{a}hler metric.

\begin{definition-new} {\rm\cite{jost,jost2}} A Hermitian manifold $(M,J,g)$ of complex dimension $m(=2n)$ is called an astheno-K\"ahler manifold if it carries a fundamental $2$-form (K\"ahler form) $\Omega$ satisfying

 $$\partial \bar{\partial} \Omega^{m-2}=0,$$ where $\partial $
 and $\bar{\partial}$ are complex exterior differentials and $\Omega^{m-2}=\underset{m-2\quad \mathrm{times}}{\Omega \wedge \ldots\wedge \Omega}$.
\end{definition-new}
The above condition is automatically satisfied for $m=2$. Thus, any Hermitian metric on a complex surface is astheno-K\"ahler. For more study \cite{Gupta}.

\subsection{Almost contact manifolds}

An odd-dimensional Riemannian manifold $\left(M^{2 n+1}, g\right)$ is said to be an almost contact metric manifold if there exists on $M$ a $(1,1)$ tensor field $\phi$, a vector field $\xi$ (called the structure vector field) and a 1-form $\eta$ such that

$$
\eta(\xi)=1, \quad \phi^2(X)=-X+\eta(X) \xi, \quad  \quad g(\phi X, \phi Y)=g(X, Y)-\eta(X)\eta(Y)
$$
for any vector fields $X, Y$ on $M$.\\
In particular, in an almost contact metric manifold, we also have $\phi \xi=0$ and $\eta \circ \phi=0$.
Such a manifold is said to be a contact metric manifold if $d \eta=\Phi$, where $\Phi(X, Y)=g(X, \phi Y)$ is called the fundamental $2 $-form of $M$. 
Such an almost contact metric structure is called a contact metric structure. Moreover, if a contact metric structure is normal, then it is called a Sasakian structure. On the other hand, if $d\Phi = 0$ and $d\eta = 0$, then $M$ is said to have an almost cosymplectic structure. In addition, if an almost cosymplectic structure is normal, then it is called a cosymplectic structure. If $N$ is a compact K\"ahler manifold, then $N \times S^1$ is the trivial example of compact cosymplectic manifolds. If a contact metric structure $(\phi, \xi, \eta)$ is normal and $d\eta = 0$, $d\Phi = 2\Phi \wedge \eta$, then $M$ is said to be Kenmotsu manifold.

These manifolds can be characterised through their Levi-Civita connection by requiring

\vspace{0.3cm}

\noindent
\hspace{3cm} $\begin{cases}
(1): \text{Sasakian } (\nabla_X \phi)Y = g(X,Y)\xi -\eta(Y)X, \\
(2): \text{Cosymplectic } (\nabla_X \phi)Y = 0, \\
(3): \text{Kenmotsu } (\nabla_X \phi)Y = g(\phi X, Y)\xi -\eta(Y)\phi X.
\end{cases}$ 

\medskip

The fundamental decomposition of the tangent bundle of an almost contact manifold $M$ splits into the direct sum of two subbundles as $TM = D \oplus \mathcal{L}$,
 where $\mathcal{L}$ is the line bundle spanned by $\xi$, which means that at any point $p\in M$, the fiber ${\mathcal L}_p$ is the $1$-dimensional subspace of $T_p M $ consisting of all scalar multiples of the vector $\xi$. This is often called the ``vertical" direction, and  $D = {\rm Ker} \eta = {\rm Im} \phi$, which means the fiber $D_p$ at a point $p$ consists of all tangent vectors $ X \in T_p M $ such that $\eta(X)=0$.  Geometrically, this is a hyperplane in the tangent space. This is often referred to as the ``horizontal" distribution or the contact distribution.

\vspace{0.5cm}
\noindent
Let $M$ be an almost contact manifold. On the product manifold $M \times \mathbb{R}$, there is a natural almost complex structure $J$ defined by
\begin{align}
J\left(X+f \frac{d}{d t}\right)=\phi X-f \xi+\eta(X) \frac{d}{d t}
\end{align}
where $X \in {\Gamma}(M)$, $t$ is the coordinate on $\mathbb{R}$, and $f$ is a smooth function on $M \times \mathbb{R}$. If $J$ is integrable, the almost contact structure is said to be normal. This is equivalent to the vanishing of the tensor field
$$
N_{\phi}:=[\phi, \phi]+d \eta \otimes \xi,
$$
where [ $\phi, \phi$ ] is the Nijenhuis torsion of $\phi$ defined by

$$
[\phi, \phi](X, Y)=[\phi X, \phi Y]+\phi^2[X, Y]-\phi[\phi X, Y]-\phi[X, \phi Y] .
$$
If the almost contact structure is normal, then
\begin{align}
\phi[\xi, X]=[\xi, \phi X] \quad \text { for all } X \in \Gamma(M) .
\label{n1}
\end{align}
In particular,

\begin{equation}
[\xi, X] \in \Gamma({\mathcal{D}}) \quad \text { for all X}  \in  {\Gamma(\mathcal{D}). }
\label{eq-1}
\end{equation}

\begin{definition-new}
{\rm \cite{oubina1}} A structure $(\phi, \xi, \eta, g)$ is a trans-Sasakian structure if and only if it is normal and

$$
d \eta=\alpha \Phi, \quad d \Phi=2 \beta \eta \wedge \Phi,
$$
where $\alpha=\dfrac{1}{2 n} \delta \phi(\xi)$, $\beta=\dfrac{1}{2 n}$ div $\xi$ and $\delta$ is the codifferential of $g$ and manifold is said to be trans-Sasakian of type $(\alpha,\beta)$.\\
A trans-Sasakian structure of type $(\alpha,\beta)$ maybe expressed as an almost contact metric structure satisfying

$$(\nabla_{X} \phi) Y=\alpha(g(X, Y) \xi-\eta(Y) X)+\beta(g(\phi X, Y) \xi-\eta(Y) \phi X ).$$
\end{definition-new}
\noindent


\noindent
It is clear that a trans-Sasakian of type $(1,0)$ is a Sasakian manifold, a trans-Sasakian of type $(0,1)$ is a Kenmotsu manifold and a trans-Sasakian of type $(0,0)$ is a cosymplectic manifold. A $(\alpha,\beta)$ trans-Sasakian manifold is an $\alpha$-Sasakian when $\alpha \in {\mathbb R},\alpha \neq 0$ and $\beta= 0$ and is $\beta$-Kenmotsu when $\beta \in {\mathbb R},\beta \neq 0$ and $\alpha= 0$.\\

\noindent
Some properties of trans-Sasakian manifolds that we will need in forthcoming sections are stated in the following lemma.
\begin{lemma}{\rm \cite{blair}}
    If $(M, \phi,\xi,\eta,g)$ is a $(\alpha,\beta)$-trans-Sasakian manifold, then
    \begin{enumerate}
        \item[\rm(i)]$ \nabla_X \xi  =-\alpha \phi X-\beta \phi^2 X,$
        \item[\rm(ii)]$\nabla_\xi X = -\alpha \phi X + \beta( X - \eta(X)\xi) + [\xi,X], \quad \nabla_\xi \xi=0,$
        \item[\rm(iii)]$\left(\nabla_X \eta\right) Y  =\alpha g(X, \phi Y)+\beta g(\phi X, \phi Y),$\\
        $(\nabla_X \phi)Y = \alpha (g(X,Y)\xi - \eta(Y)X) + \beta(g(\phi X,Y)\xi - \eta(Y)\phi X ).$
        
    \end{enumerate}
    \label{l1}
\end{lemma}
On trans-Sasakian manifold, we have
\begin{equation}
[\xi, X] \in {\mathcal D} \quad {\rm  for \ \  all}  \ X \in \Gamma (M).
\label{e1}
\end{equation}
Indeed, decomposing $X=\eta(X) \xi+X^{\mathcal{D}}$ with $X^{\mathcal{D}} \in \Gamma(\mathcal{D})$, we have

$$
[\xi, X]=\left[\xi, \eta(X) \xi+X^{\mathcal{D}}\right]=\xi(\eta(X)) \xi+\left[\xi, X^{\mathcal{D}}\right],
$$
and
$$
\xi(\eta(X))=\xi g(\xi, X)=g\left(\nabla_{\xi} \xi, X\right)+g\left(\xi, \nabla_{\xi} X\right)=-g(\xi,-\alpha \phi X-\beta \phi^2 X )=0,
$$
where we have used Lemma \ref{l1}(i) in the third and fourth equalities. Therefore

$$
[\xi, X]=\left[\xi, X^{\mathcal{D}}\right] \in \Gamma(\mathcal{D})
$$
 due to (\ref{eq-1}).

 \section{Transverse geometry of trans-Sasakian manifolds}
 Let $(M, \phi, \eta, \xi, g)$ be a trans-Sasakian manifold of dimension $2 n+1$. Recall that $\mathcal{D}=\operatorname{Ker} \eta=\operatorname{Im} \phi$ is a subbundle of $T M$ of rank $2 n$. The trans-Sasakian structure induces on $\mathcal{D}$ a natural connection $\nabla^T$, called the transverse Levi-Civita connection which, for any $U \in \Gamma(\mathcal{D})$, is defined by

\begin{align}
\nabla_{\xi}^T U=[\xi, U], \quad \nabla_X^T U=\left(\nabla_X U\right)^{\mathcal{D}}, \quad X \in \Gamma(\mathcal{D}),
\label{e2}
\end{align}
where $(\cdot)^{\mathcal{D}}$ denotes the projection onto $\mathcal{D}$. This is the only connection on $\mathcal{D}$ that satisfies

\begin{align}
\nabla_X^T\left(\left.\phi\right|_{\mathcal{D}}\right)=0, \quad \nabla_X^T\left(\left.g\right|_{\mathcal{D}}\right)=0, \quad \nabla_U^T V-\nabla_V^T U=[U, V]^{\mathcal{D}},
\label{e3}
\end{align}
for any $X \in {\Gamma}(M)$ and $U, V \in \Gamma(\mathcal{D})$. Note that

\begin{align}
\nabla_U V=& \eta(\nabla_U V) \xi+ (\nabla_U V)^{\mathcal{D}}, \quad U, V \in \Gamma(\mathcal{D}) \notag\\
=&\eta(\nabla_U V) \xi+\nabla_U^T V \notag\\
=& [-\alpha\Phi(U,V) -\beta g(U,V) + \beta\eta(U)\eta(V)]\xi + \nabla_U^T V\notag\\
=& [-\alpha\Phi(U,V) -\beta g(\phi U,\phi V)]\xi + \nabla_U^T V
\label{e4}
\end{align}
which implies

\begin{align}
[U, V]= & \eta([U, V]) \xi+[U, V]^{\mathcal{D}}, \quad U, V \in \Gamma(\mathcal{D}) \notag \\ 
=& -2\alpha \Phi(U,V)\xi+[U, V]^{\mathcal{D}}.
\label{e5}
\end{align}
Using (\ref{e2})-(\ref{e5}) we obtain the following result:

\begin{lemma}
For any $U, V, W \in \Gamma(\mathcal{D})$, the following identities hold:
\begin{enumerate}
    \item [\rm(i)] $\nabla_{[U,V]^{\mathcal{D}}}^T W = \nabla_{[U,V]}^T W + 2\alpha \Phi(U,V)[\xi,W],$
    \item[\rm(ii)] $\nabla_{[U,V]} W = 2\alpha^2\Phi(U,V)\phi W -2\alpha \beta \Phi(U,V)W  -\alpha \Phi([U,V]^{\mathcal{D}},W)\xi -\beta g([U,V],W)\xi + \nabla_{[U,V]}^T W,$
     \item[\rm(iii)] 
     
     \vspace{-.4cm}

     \begin{eqnarray}
     R(U,V)W &=& R^T(U,V)W +\alpha^2\Phi(V,W)\phi U - 2\alpha^2\Phi(U,V)\phi W -\alpha^2\Phi(U,W)\phi V \notag\\&&
     + \alpha \beta \Phi(V,W)\phi^2 U 
    +\alpha\beta g(V,W)\phi U + \beta^2 g(V,W)\phi^2 U - \alpha \beta g(U,W)\phi V \notag\\&&
    -\beta^2 g(U,W)\phi^2 V +2\alpha \beta \Phi(U,V)W - \alpha\beta \Phi (U,W)\phi^2 V,  \notag 
   \end{eqnarray}

    \item[\rm(iv)] $R (U,V)\xi = -\beta \nabla_V(\eta(V)\xi) + \beta\nabla_V(\eta(U)\xi) + \beta\eta([U,V])\xi.  $
 
\end{enumerate}
\end{lemma}
{\bf Proof.} \rm(i) is a straightforward consequence of (\ref{e5}) and (\ref{e2}): \begin{align}
        \nabla_{[U,V]}^TW = &\nabla_{(-2\alpha\Phi(U,V)\xi + [U,V]^{\mathcal{D}})}^T W \notag \\
        =& -2\alpha\Phi(U,V)[\xi,W] + \nabla_{[U,V]^{\mathcal{D}}}^T W.
    \end{align}
For \rm(ii), using (\ref{e5}) and (\ref{e4}), we compute
    \begin{align}
        \nabla_{[U,V]}W = &\nabla_{(-2\alpha\Phi(U,V)\xi + [U,V]^{\mathcal{D}})} W \notag \\
        = & -2\alpha\Phi(U,V) \nabla_\xi W + \nabla_{[U,V]^{\mathcal{D}}} W \notag \\
        = &  -2\alpha\Phi(U,V)[-\alpha\phi W +\beta (W - \eta(W)\xi) + [\xi, W]]  + \nabla_{[U,V]^{\mathcal{D}}} W \notag \\
        = &  2\alpha^2\Phi(U,V)\phi W - 2 \alpha \beta \Phi(U,V) W + 2\alpha \beta \Phi(U,V)\eta(W)\xi -2\alpha\Phi(U,V)[\xi,W] \notag\\
        &-\alpha\Phi([U,V]^{\mathcal{D}}, W)\xi -\beta g([U,V],W)\xi + \beta \eta([U,V])\eta(W)\xi + \nabla_{[U,V]}^TW \notag \\
        &+2\alpha\Phi(U,V)[\xi, W] \notag \\
        =   &  2\alpha^2\Phi(U,V)\phi W - 2 \alpha \beta \Phi(U,V) W 
        -\alpha\Phi([U,V]^{\mathcal{D}}, W)\xi\notag \\& -\beta g([U,V],W)\xi  + \nabla_{[U,V]}^TW,
\end{align}
where we have used \rm(i) in the fourth equality.\\
For \rm(iii), we use the definition $R(U,V)W = \nabla_U\nabla_V W -\nabla_V\nabla_U W -\nabla_{[U,V]} W$ and using (\ref{e2}), (\ref{e5}), Lemma \ref{l1} and \rm(ii), we obtain

    \begin{align}
R(U,V)W = &-\alpha U(\Phi(V,W))\xi - \alpha\Phi(V,W)\nabla_U\xi -\beta U (g(\phi V, \phi W))\xi - \beta g(\phi V, \phi W)\nabla_U\xi +\nabla_U^T\nabla_V^T W \notag \\
&-\alpha\Phi(U,\nabla_V^T W)\xi -\beta g(\phi U, \phi\nabla_V^T W)\xi + \alpha V(\Phi(U,W))\xi + \alpha\Phi(U,W)\nabla_V\xi + \beta V(g(\phi U,\phi W))\xi \notag \\
&+ \beta g(\phi U, \phi W)\nabla_V\xi + \alpha \Phi(V,\nabla_U^TW)\xi +\beta g(\phi V,\phi\nabla_U^T W)\xi -\nabla_V^T\nabla_U^T W  -2\alpha^2\Phi(U,V)\phi W \notag \\
&+2\alpha\beta\Phi(U,V)W + \alpha\Phi([U,V]^{\mathcal{D}},W)\xi +\beta g(\phi U,\phi W)\xi +\beta g([U,V],W)\xi - \nabla_{[U,V]}^T W \notag \\
=& R^T (U,V)W -\alpha g(\nabla_U V, \phi W)\xi -\alpha g(V,(\nabla_U\phi)W)\xi -\alpha g(V,\phi \nabla_U W)\xi -\alpha g(V,\phi W)\nabla_U \xi \notag \\
& -\beta g((\nabla_U\phi)V,\phi W)\xi -\beta g(\phi \nabla_U V,\phi W)\xi -\beta g(\phi V,(\nabla_U \phi)W)\xi -\beta g(\phi V,\phi \nabla_U W)\xi\notag \\
&-\beta g(\phi V,\phi W)\nabla_U\xi -\alpha \Phi(U,\nabla_V^T W)\xi -\beta g(\phi U,\phi\nabla_V^T W)\xi +\alpha g(\nabla_V U,\phi W)\xi \notag \\
& + \alpha g (U,(\nabla_V \phi)W)\xi +\alpha g(U, \phi\nabla_V W)\xi +\beta g((\nabla_V\phi)U,\phi W)\xi + \beta g(\phi\nabla_V U,\phi W)\xi \notag \\
&+\beta g(\phi U,(\nabla_V\phi)W)\xi + \beta g(\phi U,\phi \nabla_V W)\xi +\beta g(\phi U,phi W)\nabla_V\xi +\beta g(\phi V,\phi \nabla_U^T W)\xi \notag \\
&+\alpha\Phi(V,\nabla_U^T W)\xi + \alpha \Phi(U,W)\nabla_V \xi -2\alpha^2\Phi(U,V)\phi W + 2\alpha \beta \Phi(U,V)W +\alpha\Phi([U,V]^{\mathcal{D}},W)\xi\notag \\
&+\beta g([U,V],W)\xi \notag \\
=&R^T(U,V)W +\alpha^2\Phi(V,W)\phi U - 2\alpha^2\Phi(U,V)\phi W -
\alpha^2\Phi(U,W)\phi V  -\alpha\beta g (V,\phi W)U \notag \\
& + \alpha \beta g (V,\phi W)\eta(U)\xi +\alpha\beta g (\phi V,\phi W)\phi U -\beta^2 g (\phi V,\phi W)U + \beta^2 g(\phi V,\phi W)\eta(U)\xi \notag \\
&-\alpha\beta g(\phi U,\phi W)\phi V +\beta^2g(\phi U,\phi W)V -\beta^2 g(\phi U,\phi W)\eta(V)\xi +2\alpha\beta\Phi(U,V)W\notag \\
&+\alpha\beta\Phi(U,W)V -\alpha\beta\Phi(U,W)\eta(V)\xi \notag \\
=&R^T(U,V)W +\alpha^2\Phi(V,W)\phi U - 2\alpha^2\Phi(U,V)\phi W -\alpha^2\Phi(U,W)\phi V + \alpha \beta \Phi(V,W)\phi^2 U \notag \\
&   +\alpha\beta g(V,W)\phi U + \beta^2 g(V,W)\phi^2 U - \alpha \beta g(U,W)\phi V- \beta^2 g(U,W)\phi^2 V +2\alpha \beta \Phi(U,V)W \notag \\
& - \alpha\beta\Phi(U,W)\phi^2 V .  
 \end{align}
For \rm(iv), we compute
    \begin{align}
        R(U,V)\xi=& \nabla_U\nabla_V\xi -\nabla_V\nabla_U\xi -\nabla_{[U,V]}\xi \notag \\
        =&-\alpha \nabla_{U}\phi V +\beta \nabla_U V - \beta \nabla_U(\eta(V)\xi) +\alpha \nabla_V\phi U -\beta\nabla_V U \notag \\
        &+ \beta\nabla_V(\eta(U)\xi)+\alpha\phi[U,V]-\beta[U,V] + \beta\eta([U,V])\xi \notag \\
        =& \alpha^2\Phi(U,\phi V)\xi -\alpha\nabla_U^T\phi V -\alpha^2\Phi(V,\phi U)\xi +\alpha\nabla_V^T\phi U +\beta\nabla_U^T V -\beta\nabla_U(\eta(V)\xi) \notag \\
        &-\beta\nabla_V^T U +\beta\nabla_V(\eta(U)\xi) +\alpha\Phi[U,V] -\beta[U,V] +\beta(\eta[U,V])\xi \notag \\
        =&  -\beta \nabla_V(\eta(V)\xi) + \beta\nabla_V(\eta(U)\xi) + \beta\eta([U,V])\xi .
    \end{align}

\begin{corollary}
    For any $U,V,W \in \Gamma (\mathcal{D})$, each curvature endomorphism ${\mathcal{R}} (U,V)$ preserves $\mathcal{D}$. Moreover, for $V=\phi U$, we have that $\phi R(U,\phi U)W - R(U,\phi U)\phi W =2\alpha \beta[g(U,W)\phi U- g(U,\phi W)U]. $ 
\label{c2}
\end{corollary}

\section{Hermitian structures on the product of trans-Sasakian manifolds}

We recall next the following construction, developed independently by Tsukada \cite{tsukada} and Watson \cite{watson}, both based on a previous construction due to Morimoto \cite{morimoto}, using ideas from \cite{calabi}. With this construction, one can define a Hermitian structure on the product of two manifolds equipped with normal almost contact metric structures. We will focus later on the product of trans-Sasakian manifolds.

Let $M_1$ and $M_2$ be differentiable manifolds of dimension $2 n_1+1$ and $2 n_2+1$ and let $\left(\phi_1, \xi_1, \eta_1, g_1\right)$ and $\left(\phi_2, \xi_2, \eta_2, g_2\right)$ be almost contact metric structures on $M_1$ and $M_2$, respectively.

For $a, b \in {\mathbb{R}}, b \neq 0$, we can induce an almost Hermitian structure ( $J_{a, b}, g_{a, b}$ ) on the product manifold $M:=M_1 \times M_2$ as follows: \\

For $X_1 \in {\Gamma}\left(M_1\right)$ and $X_2 \in {\Gamma}\left(M_2\right)$, define an almost complex structure $J_{a, b}$ on $M$ \cite{tsukada} by

\begin{align}
J_{a, b}\left(X_1+X_2\right)= & \phi_1 X_1-\left(\frac{a}{b} \eta_1\left(X_1\right)+\frac{a^2+b^2}{b} \eta_2\left(X_2\right)\right) \xi_1 \notag \\
& +\phi_2 X_2+\left(\frac{1}{b} \eta_1\left(X_1\right)+\frac{a}{b} \eta_2\left(X_2\right)\right) \xi_2
\label{e6}
\end{align}
Next, define a Riemannian metric $g_{a, b}$ on $M$ by
\begin{align}
g_{a, b}\left(X_1+X_2, Y_1+Y_2\right)= & g_1\left(X_1, Y_1\right)+a\left[\eta_1\left(X_1\right) \eta_2\left(Y_2\right)+\eta_1\left(Y_1\right) \eta_2\left(X_2\right)\right] \notag  \\
& +g_2\left(X_2, Y_2\right)+\left(a^2+b^2-1\right) \eta_2\left(X_2\right) \eta_2\left(Y_2\right)
\label{e7}
\end{align}
Here, we can see that  $g_{a, b}$ is positive definite and $J_{a, b}$ is Hermitian with respect to $g_{a, b}$.

Regarding ${\Gamma}\left(M_1\right)$ and ${\Gamma}\left(M_2\right)$ as subalgebras of ${\Gamma}(M)$ in a natural manner, (\ref{e6}) and (\ref{e7}) can be rewritten in the following way, where $U_i \in \Gamma\left({\mathcal{D}}_i\right)$ :

$$
\begin{gathered}
J_{a, b} \xi_1=-\frac{a}{b} \xi_1+\frac{1}{b} \xi_2, \quad J_{a, b} U_1=\phi_1 U_1 ,\\
J_{a, b} \xi_2=-\frac{a^2+b^2}{b} \xi_1+\frac{a}{b} \xi_2, \quad J_{a, b} U_2=\phi_2 U_2,
\end{gathered}
$$
and, for $X_i, Y_i \in {\Gamma}\left(M_i\right)$ :

$$
\begin{aligned}
& g_{a, b}\left(X_1, Y_1\right)=g_1\left(X_1, Y_1\right), \quad g_{a, b}\left(X_1, X_2\right)=a \eta_1\left(X_1\right) \eta_2\left(X_2\right), \\
& g_{a, b}\left(X_2, Y_2\right)=g_2\left(X_2, Y_2\right)+\left(a^2+b^2-1\right) \eta_2\left(X_2\right) \eta_2\left(Y_2\right).
\end{aligned}
$$
Note that $g_{a, b}$ coincides with $g_1$ on $M_1$ and with $g_2$ on ${\mathcal{D}}_2$, but it modifies the length of $\xi_2$ by a factor of $a^2+b^2$; also, $\xi_1$ and $\xi_2$ are no longer orthogonal whenever $a \neq 0$. Moreover, $g_{a, b}$ is the product $g_1 + g_2$ if and only if $a=0, b= \pm 1$.

Morimoto's original construction corresponds to the case $a=0, b=1$. He proved the following result:

\begin{proposition}
    
 {\rm\cite[proposition 3]{morimoto}} Let $\left(\phi_1, \xi_1, \eta_1\right)$ and $\left(\phi_2, \xi_2, \eta_2\right)$ be almost contact structures on $M_1$ and $M_2$, respectively. Then the almost complex structure $J_{0,1}$ on $M=M_1 \times M_2$ is integrable if and only if both almost contact structures are normal.
\end{proposition}
More generally, the following result can be proved in the same way as in \cite{morimoto}:
\begin{proposition}
 Let $\left(\phi_1, \xi_1, \eta_1\right)$ and $\left(\phi_2, \xi_2, \eta_2\right)$ be almost contact structures on $M_1$ and $M_2$, respectively. If both almost contact structures are normal, then the almost complex structure $J_{a, b}$ is integrable for any $a \in {\mathbb{R}}, b \in {\mathbb{R}}, b \neq 0$.
\end{proposition}
From now on, we will deal only with the case when $( M_1, \phi_1, \xi_1, \eta_1, g_1)$  and $( M_2, \phi_2, \xi_2, \eta_2, g_2)$ are trans-Sasakian manifolds. We will denote $M_{a, b}:=M_1 \times M_2$ equipped with the Hermitian structure $\left(J_{a, b}, g_{a, b}\right)$. Moreover, we will denote simply $J:=J_{a, b},\  g:=g_{a, b}$ since there will be no risk of confusion.

In forthcoming sections, we will need explicit formulas for the Levi-Civita connection $\nabla$ on $M_{a, b}$ associated with $g$ in terms of the Levi-Civita connections $\nabla^1$ and $\nabla^2$ on $(M_1, g_1 )$ and $(M_2, g_2)$, respectively. 
Given $X_i, Y_i, Z_i \in {\Gamma}\left(M_i\right)$, we have that
\begin{align}
g\left(\nabla_{X_1} Y_1, Z_1\right)= & g_1\left(\nabla_{X_1}^1 Y_1, Z_1\right), \quad g\left(\nabla_{X_1} Y_1, Z_2\right)=a \eta_1\left(\nabla_{X_1}^1 Y_1\right) \eta_2\left(Z_2\right) \notag \\
g\left(\nabla_{X_2} Y_2, Z_1\right)= & a \eta_2\left(\nabla_{X_2}^2 Y_2\right) \eta_1\left(Z_1\right) \notag\\
g\left(\nabla_{X_2} Y_2, Z_2\right)= & g_2\left(\nabla_{X_2}^2 Y_2, Z_2\right)+\left(a^2+b^2-1\right)\left[\eta_2\left(\nabla_{X_2}^2 Y_2\right) \eta_2\left(Z_2\right)\right. \notag\\
& \left.+\eta_2\left(X_2\right) (\alpha_2g_2\left(Y_2,\phi_2 Z_2\right) + \beta_2 g_2(\phi_2 Y_2,\phi_2 Z_2))\right.\notag  \\
& \left.+\eta_2\left(Y_2\right) (\alpha_2 g_2\left( X_2,\phi_2 Z_2\right) + \beta_2 g_2(\phi_2X_2,\phi_2 Z_2))\right]\label{e8}\\
g\left(\nabla_{X_1} Y_2, Z_1\right)= &  a \eta_2\left(Y_2\right) [\alpha_1g_1\left( X_1, \phi_1 Z_1\right)+\beta_1 g_1(\phi_1 X_1,\phi_1 Z_1)]\notag\\
 g\left(\nabla_{X_1} Y_2, Z_2\right)=&a \eta_1\left(X_1\right)[\alpha_2 g_2\left(Y_2, \phi_2 Z_2\right)+\beta_2 g_2(\phi_2 Y_2, \phi_2 Z_2)] \notag\\
g\left(\nabla_{X_2} Y_1, Z_1\right)= & a \eta_2\left(X_2\right) [\alpha_1 g_1\left( Y_1, \phi_1 Z_1\right)+ \beta_1 g_1(\phi_1 Y_1,\phi_1 Z_1)]\notag \\
 g\left(\nabla_{X_2} Y_1, Z_2\right)=&a \eta_1\left(Y_1\right) [\alpha_2  g_2\left( X_2, \phi_2 Z_2\right)+ \beta_2  g_2(\phi_2 X_2,\phi_2 Z_2)]\notag 
 \label{e8}
\end{align}
The next result follows from the set of equations (\ref{e8}):
\begin{corollary} With notation as above, we have 

 $$\nabla_{X_1} Y_1=\nabla_{X_1}^1 Y_1 \in {\Gamma}\left(M_1\right),$$
$$\nabla_{X_2} Y_2=\nabla_{X_2}^2 Y_2-\left(a^2+b^2-1\right)\left[\eta_2\left(X_2\right) (\alpha_2\phi_2 Y_2+ \beta_2\phi_2^2 Y_2)+\eta_2\left(Y_2\right) (\alpha_2\phi_2 X_2+\beta_2\phi_2^2X_2)\right] \in {\Gamma}\left(M_2\right),$$
$$\nabla_{X_1} Y_2=-a\left[\alpha_1\eta_2\left(Y_2\right) \phi_1 X_1+ \alpha_2\eta_1\left(X_1\right) \phi_2 Y_2 +\beta_1\eta_2(Y_2)\phi_1^2 X_1 +\beta_2\eta_1(X_1)\phi_2^2Y_2\right] \in {\Gamma}\left(M_1\right) \oplus {\Gamma}\left(M_2\right),$$
$$\nabla_{X_2} Y_1=-a\left[\alpha_1\eta_2\left(X_2\right) \phi_1 Y_1+\alpha_2\eta_1\left(Y_1\right) \phi_2 X_2 + \beta_1\eta_2(X_2)\phi_1^2Y_1 +\beta_2\eta_1(Y_1)\phi_2^2X_2 \right] \in {\Gamma}\left(M_1\right) \oplus {\Gamma}\left(M_2\right).$$

In particular, $\nabla_{\xi_1} \xi_1=\nabla_{\xi_2} \xi_2=\nabla_{\xi_1} \xi_2=\nabla_{\xi_2} \xi_1=0$.
\label{c1}
\end{corollary}

Using the previous corollary, we compute the next $\nabla J$, which will be needed in the proof of Lemma \ref{l2} below.

\begin{lemma} For any $X_i, Y_i \in {\Gamma}\left(M_i\right), i=1,2$,
\begin{enumerate}
    \item[\rm(i)] $\left(\nabla_{X_1} J\right) Y_1=\alpha_1 g_1\left(X_1, Y_1\right) \xi_1-\alpha_1\eta_1\left(Y_1\right) X_1-\frac{a}{b} \alpha_1\Phi_1\left(X_1, Y_1\right) \xi_1+\frac{\alpha_1}{b} \Phi_1\left(X_1, Y_1\right) \xi_2\\
    -\frac{a}{b}\beta_1 g_1(\phi_1 X_1,\phi_1 Y_1)\xi_1 +\frac{\beta_1}{b}\beta_1 g_1(\phi_1 X_1,\phi_1 Y_1)\xi_2 + \beta_1 g_1(\phi_1 X_1, Y_1)\xi_1 -\beta_1\eta_1(Y_1)\phi_1 X_1\\
    -\frac{a}{b}\beta_1\eta_1(Y_1)X_1 +\frac{a}{b}\beta_1\eta_1(Y_1)\eta_1(X_1)\xi_1,$
\item[\rm(ii)] $\left(\nabla_{X_2} J\right) Y_2=\alpha_2\left[g_2\left(X_2, Y_2\right)+\left(a^2+b^2-1\right) \eta_2\left(X_2\right) \eta_2\left(Y_2\right)\right] \xi_2-\left(a^2+b^2\right)\alpha_2 \eta_2\left(Y_2\right) X_2 \\
+\beta_2 g_2\left(\phi_2X_2, Y_2\right)\xi_2 -\left(a^2+b^2\right)\beta_2 \eta_2\left(Y_2\right)\phi_2 X_2
-\frac{a^2+b^2}{b} \bigg(\alpha_2\Phi_2\left(X_2, Y_2\right)+ \beta_2 g_2(X_2,Y_2) \\-\beta_2\eta_2(X_2)\eta_2(Y_2) \bigg)\xi_1 
+\frac{a}{b}\bigg( \alpha_2\Phi_2\left(X_2, Y_2\right) + \beta_2 g_2(X_2,Y_2) -\beta_2\eta_2(X_2)\eta_2(Y_2)\bigg) \xi_2,
$

\item[\rm(iii)] $\left(\nabla_{X_1} J\right) Y_2=a \alpha_1\eta_2\left(Y_2\right) \eta_1\left(X_1\right) \xi_1-a \alpha_1\eta_2\left(Y_2\right) X_1+b\alpha_1 \eta_2\left(Y_2\right) \phi_1 X_1 -b\beta_1X_1\eta_2(Y_2)\\
+b\beta_1\eta_2(Y_2)\eta_1(X_1)\xi_1 + a\beta_1\eta_2(Y_2)\phi_1^3X_1,$
\item[\rm(iv)] $\left(\nabla_{X_2} J\right) Y_1=a\alpha_2\left[\eta_1\left(Y_1\right) \eta_2\left(X_2\right) \xi_2-\eta_1\left(Y_1\right) X_2\right]-b \alpha_2\eta_1\left(Y_1\right) \phi_2 X_2 + \frac{a^2+b^2}{b}\beta_2\phi_2^2X_2 \\+\frac{a^2}{b}\eta_1(Y_1)\beta_2\phi_2^2X_2 +a\beta_1\eta_1(Y_1)\phi_2^3X_2.$
\end{enumerate}
In particular, $\nabla_{\xi_1} J=0$ and $\nabla_{\xi_2} J=0$.
\label{l3}
\end{lemma}
{\bf Proof.} We compute $\nabla J$ using Corollary \ref{c1} and the definition of $J$.\\
For (i),
$$
\left(\nabla_{X_1} J\right) Y_1=\nabla_{X_1} J Y_1-J \nabla_{X_1}^1 Y_1.
$$
We will expand each term in detail. On the one hand,
$$
\begin{aligned}
\nabla_{X_1} J Y_1 & =\nabla_{X_1}\left(\phi_1 Y_1-\frac{a}{b} \eta_1\left(Y_1\right) \xi_1+\frac{1}{b} \eta_1\left(Y_1\right) \xi_2\right) \\
& =\nabla_{X_1}^1 \phi_1 Y_1-\frac{a}{b}\left(X_1\left(\eta_1\left(Y_1\right)\right) \xi_1-\eta_1\left(Y_1\right) \nabla_{X_1}\xi_1\right)+ \frac{1}{b}\left(X_1\left(\eta_1\left(Y_1\right)\right) \xi_2-\eta_1\left(Y_1\right) \nabla_{X_1}\xi_2\right).
\end{aligned}
$$
On the other hand,
$$ 
\begin{aligned}
    -J\nabla_{X_1}^1Y_1=-\phi_1\nabla_{X_1}^1Y_1 + \frac{a}{b}\eta_1(\nabla_{X_1}^1 Y_1)\xi_1 - \frac{1}{b}\eta_1(\nabla_{X_1}^1 Y_1)\xi_2.
\end{aligned}
$$
Putting these two expressions together, we arrive at\\
$$
\begin{aligned}
    \left(\nabla_{X_1} J\right) Y_1=&(\nabla_{X_1}^1 \phi_1) Y_1 -\frac{a}{b}g_1(\nabla_{X_1}\xi_1, Y_1)\xi_1 + \frac{1}{b}g_1(\nabla_{X_1}^1 \xi_1, Y_1)\xi_2  -\frac{a}{b}\eta_1(Y_1)\nabla_{X_1}+ \frac{1}{b}\eta_1(Y_1)\nabla_{X_1}\xi_2\\
   = &\alpha_1 g_1\left(X_1, Y_1\right) \xi_1-\alpha_1\eta_1\left(Y_1\right) X_1-\frac{a}{b} \alpha_1\Phi_1\left(X_1, Y_1\right) \xi_1+\frac{\alpha_1}{b} \Phi_1\left(X_1, Y_1\right) \xi_2 \\
   &-\frac{a}{b}\beta_1 g_1(\phi_1 X_1,\phi_1 Y_1)\xi_1 +\frac{\beta_1}{b}\beta_1 g_1(\phi_1 X_1,\phi_1 Y_1)\xi_2 + \beta_1 g_1(\phi_1 X_1, Y_1)\xi_1 -\beta_1\eta_1(Y_1)\phi_1 X_1 \\
   &-\frac{a}{b}\beta_1\eta_1(Y_1)X_1+\frac{a}{b}\beta_1\eta_1(Y_1)\eta_1(X_1)\xi_1,
\end{aligned}
$$
where we have used Lemma \ref{l1} (iii).\\
Next, for (ii)
$$
\begin{aligned}
\left(\nabla_{X_2} J\right) Y_2=&\nabla_{X_2} J Y_2-J \nabla_{X_2} Y_2\\
=&\nabla_{X_2} J Y_2 - J\bigg(\nabla_{X_2}^2 Y_2-\left(a^2+b^2-1\right)\left[\eta_2\left(X_2\right) (\alpha_2\phi_2 Y_2+ \beta_2\phi_2^2 Y_2)\right. \\
&+\left.\eta_2\left(Y_2\right) (\alpha_2\phi_2 X_2+\beta_2\phi_2^2X_2)\right] \bigg).
\end{aligned}
$$
The first term is equal to
$$
\begin{aligned}
    \nabla_{X_2} J Y_2=&\nabla(\phi_2 Y_2 -\frac{a^2+b^2}{b}\eta_2(Y_2)\xi_1+\frac{a}{b}\eta_2(Y_2)\xi_2)\\
   = &\nabla_{X_2}^2\phi_2 Y_2-(a^2+b^2-1)[\eta_2(X_2)(\alpha_2\phi_2^2Y_2 + \beta_2\phi_2^3 Y_2)+ \eta_2(\phi_2 Y_2)(\alpha_2\phi_2 X_2 + \beta_2\phi_2^2 X_2)] \\
   &-\frac{a^2+b^2}{b}\bigg(X_2(\eta_2(Y_2))\xi_1 + \eta_2(Y_2)\nabla_{X_2}\xi_1 \bigg) +\frac{a}{b}\bigg(X_2(\eta_2(Y_2))\xi_2 + \eta_2(Y_2)\nabla_{X_2}\xi_2 \bigg)\\
   =& \nabla_{X_2}^2\phi_2 Y_2 -(a^2+b^2-1)\alpha_2\eta_2(X_2)\phi_2^2Y_2 \\
   &-(a^2+b^2-1)\beta_2\eta_2(X_2)\phi_2^3Y_2 -\frac{a^2+b^2}{b}\bigg(X_2(\eta_2(Y_2))\xi_1 + \eta_2(Y_2)(\alpha_2\phi_2X_2 +\beta_2\phi_2^2X_2) \bigg) \\
   &+\frac{a}{b}\bigg(X_2(\eta_2(Y_2))\xi_2 - (a^2+b^2)\eta_2(Y_2)(\alpha_2\phi_2 X_2-\beta_2X_2 +\beta_2\eta_2(X_2)\xi_2) \bigg)
\end{aligned}
$$
and the second term is equal to
$$
\begin{aligned}
  -J\nabla_{X_2}Y_2 =&- J\bigg(\nabla_{X_2}^2 Y_2-\left(a^2+b^2-1\right)\left[\eta_2\left(X_2\right) (\alpha_2\phi_2 Y_2+ \beta_2\phi_2^2 Y_2) +\eta_2\left(Y_2\right) (\alpha_2\phi_2 X_2+\beta_2\phi_2^2X_2)\right]\bigg)\\
  =&-\phi_2(\nabla_2^2)Y_2) +(a^2+b^2-1)\bigg[\eta_2(X_2)\alpha_2\phi_2^2 Y_2 +\eta_2(X_2)\beta_2\phi_2^3 Y_2 +\alpha_2\eta_2(Y_2)\phi_2^2X_2\\
  & + \beta_2\eta_2(Y_2)\phi_2^3X_2\bigg] +\frac{a^2+b^2}{b}\eta_2(\nabla_{X_2}^2Y_2)\xi_1 -\frac{a}{b}\eta_2(\nabla_{X_2}^2Y_2)\xi_2.
\end{aligned}
$$
Therefore,
$$
\begin{aligned}
     (\nabla_{X_2}J)Y_2=& (\nabla_{X_2}^2 \phi_2)Y_2 + (a^2+b^2-1)\bigg[\alpha_2\eta_2(Y_2)\phi_2^2X_2 +\beta_2\eta_2(Y_2)\phi_2^3X_2 \bigg]+\frac{a^2+b^2}{b} g_2(\nabla_{X_2}\xi_2,Y_2)\xi_1 \\
     &+\frac{a}{b}g_2(\nabla_{X_2}\xi_2,Y_2)\xi_2\\
    =& \alpha_2\left[g_2\left(X_2, Y_2\right)+\left(a^2+b^2-1\right) \eta_2\left(X_2\right) \eta_2\left(Y_2\right)\right] \xi_2-\left(a^2+b^2\right)\alpha_2 \eta_2\left(Y_2\right) X_2 
\\
&+\beta_2 g_2\left(\phi_2X_2, Y_2\right)\xi_2-\left(a^2+b^2\right)\beta_2 \eta_2\left(Y_2\right)\phi_2 X_2
-\frac{a^2+b^2}{b} \bigg(\alpha_2\Phi_2\left(X_2, Y_2\right) + \beta_2 g_2(X_2,Y_2) \\
&-\beta_2\eta_2(X_2)\eta_2(Y_2) \bigg)\xi_1 +\frac{a}{b}\bigg( \alpha_2\Phi_2\left(X_2, Y_2\right) + \beta_2 g_2(X_2,Y_2) -\beta_2\eta_2(X_2)\eta_2(Y_2)\bigg) \xi_2
\end{aligned}
$$
using again Lemma \ref{l1} (iii).\\
Now, for (iii),
$$
\begin{aligned}
    \left(\nabla_{X_1} J\right) Y_2=&\nabla_{X_1} J Y_2-J \nabla_{X_1} Y_2\\
    =&\nabla_{X_1}(\phi_2Y_2-\frac{a^2+b^2}{b}\eta_2(Y_2)\xi_1+\frac{a}{b}\eta_2(Y_2)\xi_2) +aJ(\alpha_1\eta_2(Y_2)\phi_1X_1 +\alpha_2\eta_1(X_1)\phi_2Y_2 \\
    &+\beta_1\eta_2(Y_2)\phi_1^2X_1+\beta_2\eta_1(X_1)\phi_2^2Y_2)\\
    =& -a[\alpha_2\eta_1(X_1)\phi_2^2Y_2 + \beta_2\eta_1(X_1)\phi_2^3Y_2]-\frac{a^2+b^2}{b}[-\eta_2
    (Y_2)\alpha_1\phi_1X_1 +\beta_1X_1\eta_2(Y_2)\\
    &-\beta_1\eta_2(Y_2)\eta_1(X_1)\xi_1] -\frac{a^2}{b}\eta_2(Y_2)[\alpha_1\phi_1X_1 -\beta_1X_1 +\beta_1\eta_1(X_1)\xi_1 ] +a[\alpha_1\eta_2(Y_2)\phi_1^2X_1 \\
    &+\alpha_2\eta_1(X_1)\phi_2^2Y_2 +\beta_1\eta_2(Y_2)\phi_1^3X_1 +\beta_2\eta_1(X_1)\phi_2^3Y_2]\\
    =& a \alpha_1\eta_2\left(Y_2\right) \eta_1\left(X_1\right) \xi_1-a \alpha_1\eta_2\left(Y_2\right) X_1+b\alpha_1 \eta_2\left(Y_2\right) \phi_1 X_1 -b\beta_1X_1\eta_2(Y_2)\\
&+b\beta_1\eta_2(Y_2)\eta_1(X_1)\xi_1 + a\beta_1\eta_2(Y_2)\phi_1^3X_1.
\end{aligned}
$$
For (iv),
$$
\begin{aligned}
     \left(\nabla_{X_2} J\right) Y_1=&\nabla_{X_2} J Y_1-J \nabla_{X_2} Y_1\\
     =&\nabla_{X_2}(\phi_1Y_1-\frac{a}{b}\eta_1(Y_1)\xi_1 + \frac{1}{b}\eta_1(Y_1)\xi_2) -aJ(\alpha_1\eta_2(X_2)\phi_1Y_1 +\alpha_2\eta_1(Y_1)\phi_2X_2\\
     &+\beta_1\eta_2(X_2)\phi_1^2Y_1 +\beta_2\eta_1(Y_1)\phi_2^3X_2)\\
     =&-a\alpha_1\eta_2(X_2)\phi_1^2Y_1 -a\beta_1\eta_2(X_2)\phi_1^3Y_1 
    +\frac{a}{b}(a\eta_1(Y_1)\alpha_2\phi_2X_2 +a\eta_1(Y_1)\beta_2\phi_2^2X_2) \\
    &-\frac{1}{b}(a^2+b^2)\eta_1(Y_1)(\alpha_2\phi_2X_2 -\beta_2\phi_2^2X_2) -aJ(\alpha_1\eta_2(X_2)\phi_1Y_1 +\alpha_2\eta_1(Y_1)\phi_2X_2\\
     &+\beta_1\eta_2(X_2)\phi_1^2Y_1 +\beta_2\eta_1(Y_1)\phi_2^3X_2)\\
     =&a\alpha_2\left[\eta_1\left(Y_1\right) \eta_2\left(X_2\right) \xi_2-\eta_1\left(Y_1\right) X_2\right]-b \alpha_2\eta_1\left(Y_1\right) \phi_2 X_2 + \frac{a^2+b^2}{b}\beta_2\phi_2^2X_2 \\
     &+\frac{a^2}{b}\eta_1(Y_1)\beta_2\phi_2^2X_2 +a\beta_1\eta_1(Y_1)\phi_2^3X_2.
\end{aligned}
$$
The last statement follows easily from the previous computations.\\
In the next result, $R$ denotes the curvature tensor of $\nabla$, while $R^i$ denotes the curvature tensor of $\nabla^i, i=1,2$. We shall compute only the curvature tensors, which we will need in the proof of Theorem \ref{t1}. Let us set $\lambda_{a, b}:=a^2+b^2-1$ to shorten a little bit the statement and proof of the lemma.

\begin{lemma} With notation as above, for $U_i, V_i \in \Gamma\left({\mathcal{D}}_i\right), Z_i \in \Gamma \left(M_i\right)$,
\begin{itemize}
    \item [\rm(i)] $R\left(\xi_1, \xi_2\right)=0$,
\item[\rm(ii)]  $ R\left(U_1, V_1\right) Z_1=R^1\left(U_1, V_1\right) Z_1$ {\rm and } \\
$R\left(U_1, V_1\right) Z_2=-2 a\alpha_1\alpha_2 \Phi_1\left(U_1, V_1\right) \phi_2 Z_2 -2a\beta_2\alpha_1\Phi_1(U_1,V_1)\phi_2^2Z_2$,
\item[\rm(iii)] 
$$
\begin{aligned}
R(U_2, V_2) Z_1=&-2 a\alpha_1\alpha_2 \Phi_2(U_2, V_2) \phi_1 Z_1 -2a\alpha_2\beta_1\Phi_2(U_2,V_2)\phi_1^2Z_1 \quad {\rm  and} \\
 R(U_2, V_2) Z_2= & R^2(U_2, V_2) Z_2+\lambda_{a,b}[\Phi_2(V_2, Z_2) (\alpha_2\phi_2 U_2 + \beta_2\phi_2^2U_2)\\
 &-\Phi_2(U_2, Z_2)( \alpha_2\phi_2 V_2 +\beta_2\phi_2^2V_2) \\
&-2 \alpha_2\Phi_2(U_2, V_2) (\alpha_2\phi_2 Z_2 + \beta_2\phi_2^2Z_2)].
\end{aligned}
$$
\end{itemize}
In particular, 
$$
\begin{aligned}
R\left(U_1, V_1\right) \xi_1=-\beta_1\eta_1([U_1,V_1])\xi_1, 
R\left(U_1, V_1\right) \xi_2=0=
R\left(U_2, V_2\right) \xi_1=0,\\
R\left(U_2, V_2\right) \xi_2=R^2(U_2,V_2)\xi_2 +\lambda_{a, b}(2\alpha_2\beta_2 g_2(\phi_2U_2,\phi_2^2V_2)\xi_2 -2\beta_2^2g_2(\phi_2 U_2,\phi_2^3 V_2)\xi_2).
\end{aligned}
$$
\label{l3.}
\end{lemma}
{\bf Proof.} For (i), we compute $R\left(\xi_1, \xi_2\right) Z_i$ for $Z_i \in \Gamma \left(M_i\right), i=1,2$, using Corollary \ref{c1} and the properties of trans-Sasakian manifolds such as lemma \ref{l1}, (\ref{n1}) and (\ref{e1}). For $Z_1$, we have
$$
\begin{aligned}
R\left(\xi_1, \xi_2\right) Z_1 = &\nabla_{\xi_1} \nabla_{\xi_2} Z_1-\nabla_{\xi_2} \nabla_{\xi_1} Z_1 \\
 =&-a\nabla_{\xi_1}\left(\alpha_1 \phi_1 Z_1 +\beta_1\phi_1^2 Z_1\right)-\nabla_{\xi_2}\left(\left[\xi_1, Z_1\right]-\alpha_1\phi_1 Z_1 +\beta_1 Z_1-\beta_1\eta_1(Z_1)\xi_1\right) \\
 =&-a\left(\left[\xi_1, \alpha_1 \phi_1 Z_1 +\beta_1\phi_1^2 Z_1\right]-\varphi_1^2 Z_1\right)+a \alpha_1\phi_1\left[\xi_1, Z_1\right]- a\alpha_1\beta_1\phi_1\eta_1(Z_1)\xi_1 \\
 &-a\beta_1^2\phi_1^2\eta_1(Z_1)\xi_1 +a\beta_1\phi_1^2[\xi_1,Z_1] \\
 =&0
\end{aligned}
$$
and for $Z_2$, we have

$$
\begin{aligned}
    R\left(\xi_1, \xi_2\right) Z_2 = &\nabla_{\xi_1} \nabla_{\xi_2} Z_2-\nabla_{\xi_2} \nabla_{\xi_1} Z_2 \\
    =&\nabla_{\xi_1}(\nabla_{\xi_2}^2 Z_2-\left(a^2+b^2-1\right)\left[\eta_2\left(\xi_2\right) (\alpha_2\phi_2 Z_2+ \beta_2\phi_2^2 Z_2) +\eta_2\left(Z_2\right) (\alpha_2\phi_2 \xi_2 \right.\\
    &+\left.\beta_2\phi_2^2\xi_2)\right] )+a\nabla_{\xi_2}(\alpha_1\eta_2\left(Z_2\right) \phi_1 \xi_1+ \alpha_2\eta_1\left(\xi_1\right) \phi_2 Z_2 +\beta_1\eta_2(Z_2)\phi_1^2 \xi_1 \\
    &+\beta_2\eta_1(\xi_1)\phi_2^2Z_2)\\
    =&\nabla_{\xi_1}(\nabla_{\xi_2}^2 Z_2-\left(a^2+b^2-1\right)\eta_2\left(\xi_2\right) (\alpha_2\phi_2 Z_2+ \beta_2\phi_2^2 Z_2)  )+a\nabla_{\xi_2}( \alpha_2 \phi_2 Z_2 +\beta_2\phi_2^2Z_2)\\
    =&a\alpha_2\beta_2\phi_2\eta_2(Z_2)\xi_2 + a\beta_2^2\phi_2^2\eta_2(Z_2)\xi_2\\
    =&0.
\end{aligned}
$$Hence (i) is proved.\\
For (ii), $ R\left(U_1, V_1\right) Z_1=R^1\left(U_1, V_1\right) Z_1$ follows from the fact that $ \nabla_{X_1}Y_1=\nabla_{X_1}^1Y_1$ for $X_1, Y_1 \in \Gamma\left(M_1\right).$ For $Z_2$, we compute using Corollary \ref{c1}
$$
\begin{aligned}
   R(U_1,V_1)Z_2 = &\nabla_{U_1}\nabla_{V_1} Z_2 -\nabla_{V_1}\nabla_{U_1} Z_2 -\nabla_{[U_1,V_1]} Z_2\\
   =& -a\nabla_{U_1}(\alpha_1\eta_2(Z_2)\phi_1 V_1 + \alpha_2\eta_1(V_1)\phi_2Z_2+\beta_1\eta_2(Z_2)\phi_1^2V_1 + \beta_2\eta_1(V_1)\phi_2^2Z_2)\\
   &+a\nabla_{V_1}(\alpha_1\eta_2(Z_2)\phi_1U_1 + \alpha_2\eta_1(U_1)\phi_2Z_2 +\beta_1\eta_2(Z_2)\phi_1^2U_1 +\beta_2\eta_1(U_1)\phi_2^2Z_2)\\
   &+a(\alpha_1\eta_2(Z_2)\phi_1[U_1,V_1] + \alpha_2\eta_1([U_1,V_1])\phi_2Z_2 +\beta_1\eta_2(Z_2)\phi_1^2[U_1,V_1] \\
   & +\beta_2\eta_1([U_1,V_1])\phi_2^2Z_2)\\
   =&a\eta_2(Z_2)[-\alpha_1\nabla_{U_1}\phi_1V_1 -\beta_1\nabla_{U_1}\phi_1^2V_1 +\alpha_1\nabla_{V_1}\phi_1U_1 +\beta_1\nabla_{V_1}\phi_1^2U_1 +\alpha_1\phi_1[U_1,V_1]\\
   &\beta_1\phi_1^2[U_1,V_1] ] +a\alpha_2\eta_1([U_1,V_1])\phi_2Z_2 +a\beta_2\eta_1([U_1,V_1])\phi_2^2Z_2\\
   =& a\eta_2(Z_2)[2\alpha_1\beta_1g_1(\phi_1 U_1,\phi_1^2V_1)\xi_1 +2\beta_1\alpha_1\Phi_1(U_1,\phi_1^2V_1)\xi_1 -\alpha_1\nabla_{U_1}^T\phi_1V_1 -\beta_1\nabla_{U_1}^T\phi_1^2V_1 \\
   &+\alpha_1\nabla_{V_1}^T\phi_1U_1 +\beta_1\nabla_{V_1}^T\phi_1^2U_1 +\alpha_1\phi_1\nabla_{U_1}^{1,T}V_1 -\alpha_1\phi_1\nabla_{V_1}^{1,T}U_1 +\beta_1\phi_1^2\nabla_{U_1}^{1,T}V_1\\
   &-\beta_1\phi_1^2\nabla_{V_1}^{1,T}U_1]-2a\alpha_2\alpha_1\Phi_1(U_1,V_1)\phi_2Z_2 -2a\beta_2\alpha_1\Phi_1(U_1,V_1)\phi_2^2Z_2\\
   =&-2 a\alpha_1\alpha_2 \Phi_1\left(U_1, V_1\right) \phi_2 Z_2 -2a\beta_2\alpha_1\Phi_1(U_1,V_1)\phi_2^2Z_2
\end{aligned}
$$
due to (\ref{e3}), where $\nabla^{1,T}$ denotes the transverse connection on the trans-Sasakian manifold $M_1$.
$$
\begin{aligned}
   R(U_2,V_2)Z_1 = &\nabla_{U_2}\nabla_{V_2} Z_1 -\nabla_{V_2}\nabla_{U_2} Z_1 -\nabla_{[U_2,V_2]} Z_1\\ 
   =&-a\nabla_{U_2}[\alpha_1\eta_2(V_2)\phi_1Z_1 +\alpha_2\eta_1(Z_1)\phi_2V_2 +\beta_1\eta_2(V_2)\phi_1^2Z_1+\beta_2\eta_1(Z_1)\phi_2^2V_2]\\
   & +a\nabla_{V_2}[\alpha_1\eta_2(U_2)\phi_1Z_1 +\alpha_2\eta_1(Z_1)\phi_2 U_2 +\beta_1\eta_2(U_2)\phi_1^2Z_1 +\beta_2\eta_1(Z_1)\phi_2^2U_2]\\
   &+ a[\alpha_1\eta_2([U_2,V_2])\phi_1Z_1 +\alpha_2\eta_1(Z_1)\phi_2[U_2,V_2] + \beta_1\eta_2([U_2,V_2])\phi_1^2Z_1  \\
   & + \beta_2\eta_1(Z_1)\phi_2^2[U_2,V_2]]\\
   =& a\eta_1(Z_1)[-\alpha_2\nabla_{U_2}\phi_2V_2 -\beta_2\nabla_{U_2}\phi_2^2V_2 +\alpha_2\nabla_{U_2}\phi_2V_2 +\beta_2\phi_2^2[U_2,V_2] +\beta_2\nabla_{V_2}\phi_2^2U_2\\
   &+\alpha_2\phi_2[U_2,V_2]]+a\alpha_1\eta_2([U_2,V_2])\phi_1Z_1 +a\beta_1\eta_2([U_2,V_2])\phi_1^2Z_1\\
   =&-2 a\alpha_1\alpha_2 \Phi_2\left(U_2, V_2\right) \phi_1 Z_1 -2a\alpha_2\beta_1\Phi_2(U_2,V_2)\phi_1^2Z_1, 
\end{aligned}
$$
and 
$$
\begin{aligned}
   R(U_2,V_2)Z_2 = &\nabla_{U_2}\nabla_{V_2} Z_2 -\nabla_{V_2}\nabla_{U_2} Z_2 -\nabla_{[U_2,V_2]} Z_2\\ 
   =& \nabla_{U_2}(\nabla_{V_2}^2Z_2 -\lambda_{a,b}\eta_2(Z_2)(\alpha_2\phi_2V_2 +\beta_2\phi_2^2V_2) -\nabla_{V_2}(\nabla_{U_2}^2Z_2 -\lambda_{a,b}\eta_2(Z_2)(\alpha_2\phi_2U_2 \\
   &+\beta_2\phi_2^2U_2))-\nabla_{[U_2,V_2]}^2Z_2 +\lambda_{a,b}\eta_2([U_2,V_2])(\alpha_2\phi_2Z_2  +\beta_2\phi_2^2Z_2) \\
   &+ \lambda_{a,b}\eta_2(Z_2)(\alpha_2\phi_2[U_2,V_2] +\beta_2\phi_2^2[U_2,V_2])\\
   =& \nabla_{U_2}^2\nabla_{V_2}^2 Z_2 -\lambda_{a,b}U_2(\eta_2(Z_2))(\alpha_2\phi_2 V_2 +\beta_2\phi_2^2 V_2) -\lambda_{a,b}\eta_2(Z_2)(\alpha_2\nabla_{U_2}\phi_2 V_2 \\
   &+ \beta_2\nabla_{U_2}\phi_2^2 V_2) -\lambda_{a,b}\eta_2(\nabla_{V_2}^2 Z_2)(\alpha_2\phi_2 U_2 +\beta_2\phi_2^2 U_2) -\nabla_{V_2}^2\nabla_{U_2}^2 Z_2 \\
   &+\lambda_{a,b}\eta_2(\nabla_{U_2}^2Z_2)(\alpha_2\phi_2 V_2 +\beta_2\phi_2^2 V_2)+\lambda_{a,b} V_2(\eta_2(Z_2))(\alpha_2 \phi_2 U_2 + \beta_2 \phi_2^2 U_2) \\
   &+ \lambda_{a,b}\eta_2(Z_2)(\alpha_2\nabla_{V_2}\phi_2U_2 +\beta_2\nabla_{V_2}\phi_2^2U_2)-\nabla_{[U_2,V_2]}^2 Z_2 + \lambda_{a,b}\eta_2([U_2,V_2])(\alpha_2\phi_2 Z_2\\
   &+ \beta_2\phi_2^2 Z_2) + \lambda_{a,b}\eta_2(Z_2)(\alpha_2\phi_2[U_2,V_2]  + \beta_2\phi_2^2[U_2,V_2])\\
   =&  R^2(U_2, V_2) Z_2+\lambda_{a, b}[\Phi_2(V_2, Z_2) (\alpha_2\phi_2 U_2 + \beta_2\phi_2^2U_2)-\Phi_2(U_2, Z_2)( \alpha_2\phi_2 V_2 +\beta_2\phi_2^2V_2)\\
   & -2 \alpha_2\Phi_2(U_2, V_2) (\alpha_2\phi_2 Z_2 + \beta_2\phi_2^2Z_2)].
\end{aligned}
$$
The last statement follows from (ii), (iii) and Lemma \ref{l1} (iii).

\section{ Harmonicity of the complex structure $J_{a, b}$ with respect to $g_{a, b}$}

Let $(M, g)$ be an oriented Riemannian manifold. According to \cite{wood}, an almost complex structure $J$ on $M$ is harmonic if and only if

$$
\left[J, \nabla^* \nabla J\right]=0,
$$
where $\nabla^* \nabla J$ is the rough Laplacian of $J$ defined by $\nabla^* \nabla J=\operatorname{Tr} \nabla^2 J$. That is, if $\left\{u_1, \ldots, u_{2 n}\right\}$ is a local orthonormal frame on $M$, then

$$
\left(\nabla^* \nabla J\right)(W)=\sum_{i=1}^{2 n}\left(\nabla_{u_i, u_i}^2 J\right)(W), \quad W \in {\Gamma}(M),
$$
where the second covariant derivative of $J$ is given by

$$
\left(\nabla_{U, V}^2 J\right)(W)=\left(\nabla_U\left(\nabla_V J\right)\right)(W)-\left(\nabla_{\nabla_U V} J\right)(W).
$$
It is clear that $\nabla^* \nabla J$ is a $(1,1)$-tensor on $M$.\\
Let ( $M^{2n}, J, g$ ) be an almost Hermitian manifold. According to \cite{wood}, the following 2-form $\rho$ plays a special role in determining whether the almost complex structure $J$ is harmonic:

$$
\rho=R(\omega) \in \Omega^2(M),
$$
where $\omega$ is the fundamental 2-form associated with ( $J, g$ ) and $R$ is the curvature operator acting on 2-forms. This 2-form $\rho$ is a natural generalization of the Chern-Ricci form of a K\"ahler manifold, although in general it is not closed. It can be seen that the skew-symmetric tensor $P: T M \rightarrow T M$ obtained by contracting $\rho$ and $g$, that is, $\rho(X, Y)=g(P X, Y)$, is given by

\begin{equation}
PX=\frac{1}{2} \sum_{i=1}^{2 n} R\left(e_i, J e_i\right) X, \quad X \in {\Gamma}(M),
\label{n2}
\end{equation}
where $\left\{e_i\right\}$ is any orthonormal local frame of $M$. \\

Let $\delta J \in {\Gamma}(M)$ denote the codifferential of $J$, that is, the unique vector field on $M$ satisfying

$$
g(\delta J, X)=\delta \omega(X) \quad \text { for all } X \in {\Gamma}(M),
$$
where $\delta \omega$ is the codifferential of $\omega$. Since $\delta \omega$ is given by

$$
\delta \omega(X)=-\sum_{i=1}^{2 n}\left(\nabla_{e_i} \omega\right)\left(e_i, X\right)
$$
for any local orthonormal frame $\left\{e_i\right\}$ of $M$, we obtain the following expression for $\delta J$ :

\begin{align}
\delta J=\sum_{i=1}^{2 n}\left(\nabla_{e_i} J\right)\left(e_i\right).
\label{n3}
\end{align}
With all these ingredients, we may recall the following result from \cite{wood}.

\begin{proposition}
 {\rm \cite [Theorem 2.8]{wood}} Let $J$ be the almost complex structure of a $2 n$-dimensional almost Hermitian manifold $(M, J, g)$. If $J$ is integrable, then

$$
\left[J, \nabla^* \nabla J\right]=2\left(\nabla_{\delta J} J-[J, P]\right)
$$

In particular, $J$ is harmonic if and only if $[J, P]=\nabla_{\delta J} J$.
\label{p1}
\end{proposition}

First, we prove an auxiliary result, which generalizes \cite[Lemma 5.4]{tsukada}, where the case of Calabi-Eckmann manifolds is considered.

\begin{lemma}Let $(M_1^{2 n_1+1}, \phi_1, \xi_1, \eta_1, g_1)$ and $( M_2^{2 n_2+1}, \phi_2, \xi_2, \eta_2, g_2 )$  be two trans-sasakian manifolds. If $(J, g):=\left(J_{a, b}, g_{a, b}\right)$ denotes the complex structure on $M_{a, b}=M_1 \times M_2$ given in {\rm(\ref{e6})} and {\rm(\ref{e7})}, then the codifferential $\delta J$ of the complex structure $J$ on $M_{a, b}$ is given by

$$
\delta J=2 n_1( \alpha_1\xi_1 -\frac{a}{b}\beta_1\xi_1 +\frac{\beta_1}{b}\xi_2)+2 n_2 (\alpha_2\xi_2 +\beta_2\xi_1 +\frac{a}{b}\beta_2\xi_2).
$$
Moreover, $\nabla_{\delta J} J=0$.
\label{l2}
\end{lemma}
{\bf Proof.} We will compute $\delta J$ using (\ref{n3}). Let us consider a local orthonormal frame on $M_{a, b}$ of the following form:

$$
\left\{\xi_1, J \xi_1=-\frac{a}{b} \xi_1+\frac{1}{b} \xi_2, e_1, \ldots, e_{2 n_1}, f_1, \ldots, f_{2 n_2}\right\},
$$
where each $e_j$, $j=1,\ldots, n_1$ is a local section of ${\mathcal{D}}_1$ and each $f_k$, $k=1,\ldots, n_2$ is a local section of ${\mathcal{D}}_2$. With this frame, (\ref{n3}) becomes

$$
\delta J=\left(\nabla_{\xi_1} J\right) \xi_1+\left(\nabla_{J \xi_1} J\right) J \xi_1+\sum_{j=1}^{2 n_1}\left(\nabla_{e_j} J\right) e_j+\sum_{k=1}^{2 n_2}\left(\nabla_{f_k} J\right) f_k.
$$
Since $J$ is integrable, it follows from \cite[ Corollary 2.2]{gray} that $\left(\nabla_{J \xi_1} J\right) J \xi_1=\left(\nabla_{\xi_1} J\right) \xi_1$. Therefore, it follows from Lemma \ref{l3} that $$\left(\nabla_{\xi_1} J\right) \xi_1=\left(\nabla_{J \xi_1} J\right) J \xi_1=0,$$ $$\left(\nabla_{e_j} J\right) e_j=\alpha_1\xi_1 -\frac{a}{b}\beta_1\xi_1 +\frac{\beta_1}{b}\xi_2$$ and $$\left(\nabla_{f_k} J\right) f_k=\alpha_2\xi_2 +\beta_2\xi_1 +\frac{a}{b}\beta_2\xi_2,$$ 
so that
$$
\begin{aligned}
\delta J=&\sum_{j=1}^{2 n_1}\alpha_1\xi_1 -\frac{a}{b}\beta_1\xi_1 +\frac{\beta_1}{b}\xi_2 +\sum_{k=1}^{2 n_2} \alpha_2\xi_2 +\beta_2\xi_1 +\frac{a}{b}\beta_2\xi_2\\
=&2 n_1( \alpha_1\xi_1 -\frac{a}{b}\beta_1\xi_1 +\frac{\beta_1}{b}\xi_2)+2 n_2 (\alpha_2\xi_2 +\beta_2\xi_1 +\frac{a}{b}\beta_2\xi_2)
\end{aligned}
$$
as stated. Finally, the last statement follows from Lemma \ref{l3}.
\\

The main result of this section is the following:
\begin{theorem}
 Let $(M_1^{2 n_1+1}, \phi_1, \xi_1, \eta_1, g_1)$ and $( M_2^{2 n_2+2}, \phi_2, \xi_2, \eta_2, g_2 )$  be two trans-sasakian manifolds. If $(J, g):=\left(J_{a, b}, g_{a, b}\right)$ denotes the complex structure on $M_{a, b}=M_1 \times M_2$ given in {\rm(\ref{e6})} and {\rm(\ref{e7})}, then $J$ is harmonic with respect to $g$ 
 if $$2\alpha_1\beta_1[g_1(e_j,U_1)\phi_1e_j -g_1(e_j,\phi_1U_1)e_j]=0$$
 and 
 $$2\alpha_2\beta_2[g_2(f_k,U_2)\phi_2f_k-g_2(f_k,\phi_2U_2)f_k]=0$$
 for any $a, b \in {\mathbb{R}}, b \neq 0$.
 \label{t1}
\end{theorem}

\noindent{\textbf{Proof}.} As in Lemma \ref{l2}, we consider a local orthonormal frame on $M_{a, b}$ of the following form:

\begin{align}
\left\{\xi_1, J \xi_1=-\frac{a}{b} \xi_1+\frac{1}{b} \xi_2, e_1, \ldots, e_{2 n_1}, f_1, \ldots, f_{2 n_2}\right\}.
\label{n4}
\end{align}
Since $J$ is integrable, it follows from Proposition \ref{p1} and Lemma \ref{l2} that $J$ is harmonic if and only if $[J, P]=0$, with $P$ defined as in (\ref{n2}) for the local frame (\ref{n4}). We prove next that $J$ and $P$ commute indeed. We begin by computing $P$ :

$$
\begin{aligned}
2 P & =2 R\left(\xi_1, J \xi_1\right)+\sum_j R\left(e_j, \phi_1 e_j\right)+\sum_k R\left(f_k, \phi_2 f_k\right) \\
& =\sum_j R\left(e_j, \phi_1 e_j\right)+\sum_k R\left(f_k, \phi_2 f_k\right)
\end{aligned}
$$
due to Lemma \ref{l3} .\\
Recalling that $R\left(e_j, \phi_1 e_j\right) \xi_i= R\left(f_k, \phi_1 f_k\right) \xi_i=0, \  i=1,2$.
We use Lemma \ref{l3} for all the computations below. First, for $U_1 \in \Gamma\left({\mathcal{D}}_1\right)$, we compute

$$
\begin{aligned}
\left[J, R\left(e_j, \phi_1 e_j\right)\right]\left(U_1\right)=&J\left(R^1\left(e_j, \varphi_1 e_j\right) U_1\right)-R^1\left(e_j, \phi_1 e_j\right) J U_1\\
=&\phi_1(R^1(e_j,\phi_1e_j)U_1) -R^1(e_j,\phi_1e_j)\phi_1U_1\\
=&2\alpha_1\beta_1[g_1(e_j,U_1)\phi_1e_j -g_1(e_j,\phi_1U_1)e_j]
\end{aligned}
$$
due to Corollary \ref{c2}. Now, for $U_2 \in \Gamma\left({\mathcal{D}}_2\right)$,
$$
\begin{aligned}
{\left[J, R\left(e_j, \phi_1 e_j\right)\right]\left(U_2\right) }  =&J(-2 a \alpha_1\alpha_2\Phi_1\left(e_j, \phi_1 e_j\right) \phi_2U_2 -2a\beta_2\alpha_1\Phi_1(e_j,\phi_1e_j)\phi_2^2U_2)\\
&-R\left(e_j, \phi_1 e_j\right) \phi_2 U_2 \\
=&-2 a \alpha_1\alpha_2\Phi_1\left(e_j, \phi_1 e_j\right) \phi_2^2U_2 -2a\beta_2\alpha_1\Phi_1(e_j,\phi_1e_j)\phi_2^3U_2-R\left(e_j, \phi_1 e_j\right) \phi_2 U_2\\
=&0.
\end{aligned}
$$
Similarly
$$
\begin{aligned}
    {\left[J, R\left(f_k, \phi_2 f_k\right)\right]\left(U_1\right) }=&J
    (R\left(f_k, \phi_2 f_k\right)U_1)-R\left(f_k, \phi_2 f_k\right)\phi_1U_1\\
    =&J(-2a\alpha_1\alpha_2\Phi_2(f_k,\phi_2f_k)\phi_1U_1 -2a\alpha_2\beta_1\Phi_2(f_k,\phi_2 f_k)\phi_1^2U_1)\\
    &-(-2a\alpha_1\alpha_2\Phi_2(f_k,\phi_2 f_k)\phi_1^2U_1 -2a\alpha_2\beta_1\Phi_2(f_k,\phi_2 f_k)\phi_1^3U_1)\\
    =&0.
\end{aligned}
$$
Finally,
$$
\begin{aligned}
    {\left[J, R\left(f_k, \phi_2 f_k\right)\right]\left(U_2\right) }=&J
    (R\left(f_k, \phi_2 f_k\right)U_2)-R\left(f_k, \phi_2 f_k\right)\phi_2U_2\\
    =&J\bigg(R^2(f_k,\phi_2f_k)U_2 + \lambda_{a,b}[\Phi_2(\phi_2f_k,U_2)(\alpha_2\phi_2f_k +\beta_2\phi_2^2f_k)\\
    &-\Phi_2(f_k,U_2)(\alpha_2\phi_2^2f_k +\beta_2\phi_2^3f_k) -2\alpha_2\Phi_2(f_k,\phi_2f_k)(\alpha_2\phi_2U_2 +\beta_2\phi_2^2U_2)])\\
    &-\lambda_{a,b}\eta_2(U_2)[-2\alpha_2\beta_2 g_2(\phi_2f_k,\phi_2^3f_k)\xi_2 + 2\beta_2^2g_2(\phi_2f_k,\phi_2^4f_k)\xi_2]\bigg)\\
    =&J\bigg(R^2(f_k,\phi_2f_k)\phi_2U_2 + \lambda_{a,b}[\Phi_2(\phi_2f_k,\phi_2U_2)(\alpha_2\phi_2f_k +\beta_2\phi_2^2f_k)\\
    &-\Phi_2(f_k,\phi_2U_2)(\alpha_2\phi_2^2f_k +\beta_2\phi_2^3f_k) -2\alpha_2\Phi_2(f_k,\phi_2f_k)(\alpha_2\phi_2^2U_2 +\beta_2\phi_2^3U_2)])\\
    &-\lambda_{a,b}\eta_2(\phi_2U_2)[-2\alpha_2\beta_2 g_2(\phi_2f_k,\phi_2^3f_k)\xi_2 + 2\beta_2^2g_2(\phi_2f_k,\phi_2^4f_k)\xi_2]\bigg)\\
    =&\phi_2(R^2(f_k,\phi_2f_k)U_2)- R^2(f_k,\phi_2f_k)\phi_2U_2\\
    =&2\alpha_2\beta_2[g_2(f_k,U_2)\phi_2f_k-g_2(f_k,\phi_2U_2)f_k]
\end{aligned}
$$
due to Corollary \ref{c2}.\\
\vspace{.5cm}
\noindent
\begin{tabular}{|p{.7 cm}|p{2.5cm}|p{2.5cm}|p{.3 cm}|p{.3 cm}|p{.3 cm}|p{.3 cm}|p{3cm}| }
 \hline
 \multicolumn{8}{|c|}{Table 1: Particular cases} \\
 \hline
 S.No. &$M_1$ & $M_2$& $\alpha_1$ & $\alpha_2$ & $\beta_1$&$\beta_2$& 
 Harmonicity \\
 \hline
1&$\alpha_1$-Sasakian & $\alpha_2$-Sasakian  & 1 & 1 &0  & 0 & Yes \\
\hline
2&$\alpha_1$-Sasakian & $\beta_2$-Kenmotsu  & 1 &0  & 0 &1  &Yes \\
\hline
3&$\alpha_1$-Sasakian & Cosymplectic   &1  & 0 & 0 &0  &Yes   \\
\hline
4&$\beta_1$-Kenmotsu & $\beta_2$-Kenmotsu  & 0 & 0 & 1 &1  & 
Yes  \\
\hline
5&$\beta_1$-Kenmotsu & $\alpha_2$-Sasakian &0  & 1  & 1 & 0 & 
Yes   \\
\hline
6&$\beta_1$-Kenmotsu & Cosymplectic & 0 & 0 & 1 & 0 & Yes   \\
\hline
7&Cosymplectic & $\alpha_2$-Sasakian  & 0 & 1 & 0 & 0 & Yes \\
\hline
8&Cosymplectic &  $\beta_2$-Kenmotsu  & 0 & 0 & 0 & 1 & Yes   \\
\hline
9&Cosymplectic &   Cosymplectic  & 0 & 0 & 0 & 0 & Yes  \\
 \hline
 \end{tabular}\\
\vspace{.3cm}
\noindent
Since $\left[J, R\left(e_j, \varphi_1 e_j\right)\right]=\left[J, R\left(f_k, \varphi_2 f_k\right)\right]=0$ in all the above cases. It follows that $[J,P]=0$ and thus the complex structure $J=J_{a,b}$ is harmonic with respect to the metric $g=g_{a,b}$.

\begin{remark}
    For case {\rm 1} of the above table when  $\alpha_i=1$ for $i=1,2$, it is proved in {\rm \cite{andrada}} that product manifold have harmonic complex structure with respect to its metric.
\end{remark}

\begin{theorem}
   {\rm\cite{marrero}}  A $(\alpha, \beta)$-trans-Sasakian manifold of dimension
 $\geq 5$  is either $\alpha$-Sasakian, $\beta$-Kenmotsu or cosymplectic.
 \label{t2}
\end{theorem}

\begin{result}
 In case $\rm 9$ of the above table, from {\rm \cite{capursi}}, the product of two cosymplectic manifolds will always be an astheno-K\"ahler manifold and from theorem \ref{t1} we can see its complex structure is also harmonic with respect to its metric.
\end{result}

\begin{result}
    Again, for case {\rm 3} and {\rm 4} of the above table, we can observe from {\rm \cite{nidhi}} that if $\dim_{\mathbb{C}} M_{a,b}=m\geq 3$, a product of $\alpha$-Sasakian and cosymplectic is astheno-K\"ahler only when the  dimension of $\alpha$-Sasakian is $3$ so from theorem \ref{t1},  complex structure on this product manifold is also harmonic with respect to its metric.
\end{result}

\begin{result}
   In case {\rm 1} of the above table, product of two $\alpha$-Sasakian manifolds is also astheno-K\"ahler when 
   $\dim_{\mathbb{C}} M_{a,b}= m=3$ and from Theorem \ref{t1}, the complex structure on it is harmonic with respect to its metric as well.
\end{result}
\noindent
\textbf{Conclusion:} In this study, we have concluded that in many of the cases astheno-K\"ahler \cite{nidhi} have a harmonic complex structure with respect to its metric as well and from theorem \ref{t2} the product of two $(\alpha, \beta)-$trans-Sasakian manifold $M_{a, b}=M_1 \times M_2$ with complex structure $(J, g):=\left(J_{a, b}, g_{a, b}\right)$ particularly reduced to the cases of table 1 when $\dim_{\mathbb{R}} M_1 \geq 5$ and $\dim_{\mathbb{R}} M_2 \geq 5$  eventually have harmonic complex structure.

 
\medskip 
\noindent \textbf {Declarations:}\\
\textbf { Data Availability:} No new data were created or analysed in this study.\\
\textbf {Conflicts of interest:} The authors declare that they have no known financial conflicts of
interest or personal relationships that could have influenced the work presented in this paper.\\
\textbf {Funding:} The corresponding author is supported and funded by the National Board of Higher Mathematics (NBHM) project no. 02011/21/2023NBHM(R.P.)/R\&DII/14960, INDIA.\\
\textbf {Ethics approval:} The authors hereby affirm that the contents of this manuscript are original.
Furthermore, it has neither been published elsewhere in any language fully or partly, nor is it under
review for publication anywhere.

\end{document}